\documentclass[preprint, ECP]{ejpecp} 

\usepackage{xcolor}

\SHORTTITLE{How smooth is the drift of the mixed fBm?}

\TITLE{How smooth is the drift of the mixed fractional Brownian motion?\support{supported
    by the ISF grant 2821/25}\ }

\AUTHORS{%
  Pavel Chigansky\footnote{The Hebrew University of Jerusalem, Israel
    \EMAIL{Pavel.Chigansky@mail.huji.ac.il}}\orcid{0000-0002-3532-5783}
  \and 
  Marina Kleptsyna\footnote{Le Mans Universit\'{e}, France. 
  \EMAIL{Marina.Kleptsyna@univ-lemans.fr}}\orcid{0000-0002-3459-8413}} 


\AMSSUBJ{60G44; 60G22} 

\SUBMITTED{December 3, 2025} 
\ACCEPTED{March 27, 2026} 

\ARXIVID{2511.22542} 

\VOLUME{0}
\YEAR{2026}
\PAPERNUM{0}
\DOI{10.1214/YY-TN}

\ABSTRACT{
The mixed fractional Brownian motion - the sum of independent fractional and standard Brownian motions - is known to be a semimartingale if the Hurst exponent $H$ of its fractional component satisfies $H > 3/4$. The question posed in the title is motivated by recent findings in quantitative finance. In this note, we show that the drift in its Doob–Meyer decomposition has a derivative that is $\gamma$-Hölder continuous for any $\gamma < 2H - 3/2$. }

\newcommand{\Real}{\mathbb R}
\newcommand{\cov}{\mathsf{Cov}}
\newcommand{\F}{\mathcal{F}}
\newcommand{\E}{\mathsf{E}}

\begin{document}



\section{The main result}
The fractional Brownian motion (fBm) $B^H=(B^H_t, t\in \Real_+)$ is the Gaussian process with continuous paths, zero mean and covariance function 
$$
\cov(B^H_s,B^H_t) = \frac 1 2 \big(t^{2H} +    s^{2H} -  |t-s|^{2H}\big),
$$
where $H\in (0,1)$ is its Hurst exponent (\cite{EM02}, \cite{M08}, \cite{PT17}).
While $B^H$ is not a semimartingale for $H\ne 1/2$, it can become one if an independent Brownian motion $B=(B_t, t\in \Real_+)$ is added to it.
More precisely, the process
\begin{equation}\label{mfBm}
X_t = B^H_t + B_t, \quad t\in \Real_+,
\end{equation}
called  mixed fBm, is a semimartingale if and only if $H\in \{1/2\}\cup (3/4, 1]$.
This effect, discovered in \cite{Ch01}, plays a role in financial applications, \cite{Ch03}, \cite{BSV07}.
Some other interesting features of this process were revealed in, e.g.,  \cite{BN03},  \cite{CCK}, \cite{ChKM20}, \cite{CDM2025}, \cite{DSSV24},  \cite{vZ07}.

By the Doob-Meyer theorem, the continuous semimartingale $X$ is uniquely decomposable into the sum of a local martingale and a continuous process of locally finite variation, referred to in this note as the drift. The analysis of the drift of the mixed fBm is motivated by recent findings in quantitative finance. Two key concepts in this field are volatility, representing the intensity of price fluctuations, and market impact, defined as the magnitude of the average price return induced by a given volume of transactions; see \cite{gatheral2022volatility, jusselin2020no} for details. It was recently demonstrated in \cite{muhle2025metaflow} that the mixed fBm with $H > 3/4$ serves as a natural model for the signed flow of financial transactions.   As it turns out, the smoothness of the drift determines the exact behavior of the price volatility; see \cite[Theorem 4.1]{muhle2025metaflow}.

\medskip
\noindent
The main result of this note is the precise value of its H\"older exponent.

\begin{theorem}\label{thm}
The drift of the mixed fBm with $H>3/4$ is differentiable and the derivative has a continuous modification which, with probability one, 
is locally $\gamma$-H\"older continuous for all $\gamma<2H-3/2$.   No such modification exists for $\gamma > 2H-3/2$. 
\end{theorem}

\begin{remark}
The question of existence of the $\gamma$-H\"older continuous modification for $\gamma=2H-3/2$ requires a more delicate analysis and is left for further research.
\end{remark}
 
\section{The proof of Theorem \ref{thm}}

\subsection{The drift process.}
The semimartingale decomposition of the mixed processes such as \eqref{mfBm} was obtained in \cite{Hi68}, \cite{Sh66}. A comprehensive discussion of this theory in the context of mixed fBm can be found in \cite{Ch03b}.
For $H\in (3/4, 1]$, the mixed fBm, defined in \eqref{mfBm}, is a semimartingale in its own filtration $\F^X_t = \sigma\{X_s, s\in [0,t]\}$
with the Doob-Meyer decomposition 
\begin{equation}\label{DM}
X_t = \overline B_t - \int_0^t \varphi_s(X) ds, \quad t\in \Real_+,
\end{equation}
where $\overline B$ is a Brownian motion adapted to $\F^X$. The drift process is given by the stochastic integral  
\begin{equation}\label{phi}
\varphi_s (X) = \int_0^s L(r,s)dX_r, 
\end{equation}
where the kernel $L(r,s)$ is the unique $L^2([0,s])$ solution to the integral equation 
\begin{equation}\label{Leq}
L(r,s) +  \int_0^s L(\tau,s) c_H|r-\tau|^{2H-2}d\tau = -c_H |s-r|^{2H-2}, \quad 0<r<s,
\end{equation}
with the constant $c_H =  H(2H-1)$.

\begin{remark}

The assertion of Theorem \ref{thm} is obvious for $H=1$, since in this case  \eqref{Leq} reduces to the explicitly solvable equation
$$
L(r,s) +  \int_0^s L(\tau,s)  d\tau = -1 , \quad 0<r<s.
$$ 
Its unique solution is constant in $r$ and is given by $L(r,s)=-1/(1+s)$.  
Consequently, the process in \eqref{phi} takes the  form
$\varphi_s = -X_s/(1+s)$. 
It has the same H\"older exponent as the mixed fBm itself, which, for $H=1$,
reduces to
\footnote{
This is not to be confused with the representation \eqref{DM} in which $\overline B$ is a Brownian motion adapted to the filtration generated by $X$.
}
 $X_s= \xi s+ B_s$ where the random variable $\xi\sim N(0,1)$ is independent of $B$.
Thus, $X$ is as smooth as the Brownian motion - that is, H\"older continuous for any $\gamma < 1/2$, as claimed. 
The case $H\in (3/4, 1)$ is more subtle, since \eqref{Leq} is no longer solvable in a closed form.
\end{remark}

\begin{remark}
A rough  intuition for the Hölder continuity of $\varphi_t(X)$ is as follows. The second term on the left-hand side of equation \eqref{Leq} is smoother than the first due to the integration; therefore, its solution $L(r,s)$ is as smooth as the forcing function on the right-hand side, possessing a square-integrable singularity at the endpoint $s$. For $H > 3/4$, the fractional component $B^H$ of $X$ is smoother than its  Brownian component $B$, rendering $X$  $\gamma$-Hölder continuous for any $\gamma<1/2$. Just as the Riemann integration in 
$$
\int_0^t (t-s)^{-\alpha}ds=\frac 1{1-\alpha} t^{1-\alpha}, \qquad \alpha\in (0,1),
$$
increases the smoothness of the integrand $(t-s)^{-\alpha}$ by $1$, the stochastic integration in \eqref{phi} increases the smoothness of $L(r,s)$ by $1/2$. 
Consequently, the resulting process $\varphi_t(X)$ can be expected to have a Hölder exponent arbitrarily close to 
$2H - 2 + 1/2 = 2H - 3/2$. 
\end{remark}

To prove the claimed H\"older continuity we show that, for $H\in (3/4, 1)$,
\begin{equation}\label{proveme}
\E (\varphi_t(X)-\varphi_s(X))^2\le C_2(s,t) |t-s|^{4H-3}
\end{equation}
for some bounded function $C_2(s,t)$.
Since the process $\varphi(X)$ in \eqref{phi} is Gaussian, this implies  
$$
\E |\varphi_t(X)-\varphi_s(X)|^{2m}\le C_m(s,t)|t-s|^{(4H-3)m}, \quad \forall m\ge 1,
$$
with bounded functions $C_m(s,t)$. Then by Kolmogorov's continuity theorem, $\phi(X)$ has a continuous modification, which is a.s. locally $\gamma$-H\"older
continuous for all 
$$
\gamma < \frac{m(4H-3)-1}{2m}.
$$
The assertion of Theorem \ref{thm} follows by arbitrariness of $m$.

To prove that $\gamma$-H\"older continuous modification does not exist for $\gamma > 2H-3/2$ we will show that 
for some $t>0$, 
\begin{equation}\label{lower}
\varliminf_{s\to t}\frac{1}{|t-s|^{4H-3}}\E (\varphi_t(X)-\varphi_s(X))^2 >0.
\end{equation}
We can write  
$$
\frac{\varphi_t(X)-\varphi_s(X)}{\sqrt{\E (\varphi_t(X)-\varphi_s(X))^2}} =  
\frac{\varphi_t(X)-\varphi_s(X)}{|t-s|^{\gamma}}\sqrt{ \frac{|t-s|^{4H-3}}{\E (\varphi_t(X)-\varphi_s(X))^2 }  }|t-s|^{\gamma - (2H-3/2)}.
$$
If $\phi(X)$ is a $\gamma$-H\"older continuous modification with $\gamma> 2H-3/2$ and \eqref{lower} holds, the right hand side converges a.s. to zero as $s\to t$. 
However, the random variable in the left hand side has the standard normal distribution for all $s\ne t$. The obtained contradiction verifies the claim.

\subsection{Auxiliary estimates}
To prove \eqref{proveme}, we will need some estimates for the solutions to weakly singular equations such as \eqref{Leq}.
For brevity, define 
$$
\alpha =2-2H\in (0,1/2), \qquad b_\alpha:= c_{1-\alpha/2}=(1-\alpha/2)(1-\alpha).
$$
For a fixed $t>0$, let $f(t;r,s)$, $r,s\in [0,t]^2$ be a real valued function with finite $L^2([0,s])$ norm,
$$
\|f(t;\cdot,s)\|_2^2 = \int_0^s f(t;r,s)^2 dr<\infty, \quad \forall s\in [0,t].
$$
Consider the integral equation 
\begin{equation}\label{Qeq}
Q(r,s) +  \int_0^s Q(\tau,s) b_\alpha|r-\tau|^{-\alpha} d\tau = f(t; r,s), \quad 0< r< s< t.
\end{equation}
For each fixed $s$ and $t$, it has the unique Hilbert-Schmidt solution $Q(\cdot, s)\in L^2([0,s])$, see \cite[Section 94]{RN55}, since, for $\alpha \in (0,1/2)$, the kernel satisfies
$$
\int_0^s\int_0^s |r-\tau|^{-2\alpha}drd\tau<\infty.
$$ 
Since the integral in \eqref{Qeq} depends only on the values of $Q(\tau,s)$ for $\tau<s$, this equation can be used to extend the definition of the function $Q(r,s)$ to $r>s$. This extension satisfies \eqref{Qeq} for all $r,s\in [0,t]$.

The following lemma establishes that the solution to \eqref{Qeq} inherits from the forcing function $f$ the types of regularity relevant to our analysis.

\begin{lemma}\label{lem}
\

\medskip

\begin{enumerate}

\renewcommand{\theenumi}{\roman{enumi}}
\item\label{i} Assume that $\|f(t;\cdot,s)\|_\infty := \sup_{r\in [0,t]}|f(t;r,s)|<\infty$, then there exists a locally bounded function $C_1(t)$ such that
$$
 |Q(s,r)| \le C_1(t)\|f(t;\cdot,s)\|_\infty, \quad \forall r,s \in [0,t].
$$

\item\label{ii} Assume that $|f(t;r,s)|\le C(s,t) |s-r|^{-\alpha}$ for some function $C(s,t)$ which satisfies  
$$
\sup_{s\in [0,t]}C(s,t)<\infty.
$$
Then there exists a locally bounded function  $C_2(t)$ such that
\begin{equation}\label{Qstalpha}
|Q(r,s)| \le     C_2(t) C(s,t)  |s-r|^{-\alpha}, \quad \forall r,s \in [0,t].
\end{equation}

\item\label{iii} Let $f(t;r,s)= (s-r)^{-\alpha}-(t-r)^{-\alpha}$, then there exists a locally bounded function $C_3(t)$ such that
$$
|Q(r,s)|\le C_3(t)\Big((t-s)^{1-2\alpha}+(s-r)^{-\alpha}-(t-r)^{-\alpha}\Big), \qquad  0<r<s<t.
$$
\end{enumerate}

\end{lemma}

\begin{proof}\

\medskip
\noindent
\eqref{i} Multiply both sides of \eqref{Qeq} by $Q(r,s)$ and integrate to obtain 
$$
\big\|Q(\cdot,s)\big\|^2_2 +   \int_0^s\int_0^s Q(r,s) Q(\tau,s) b_\alpha|r-\tau|^{-\alpha} d\tau dr =\int_0^s Q(r,s) f(t; r,s)dr.
$$
The kernel $|r-\tau|^{-\alpha}$ is non-negative definite and hence
$$
\|Q(\cdot,s)\|^2_2\le \int_0^s |Q(r,s)| |f(t; r,s)| dr \le \|f(t;\cdot, s)\|_\infty \|Q(\cdot,s)\|_2\sqrt{t}.
$$
 Since $Q(\cdot,s)\in L^2([0,s])$ it follows that 
$\|Q(\cdot,s)\|_2 \le \|f(t;\cdot, s)\|_\infty  \sqrt{t}$ 
and, consequently, it follows from \eqref{Qeq} that  
\begin{align*}
 |Q(r,s)| \le\ & \|f(t;\cdot, s)\|_\infty  + 
b_\alpha \|Q(\cdot,s)\|_2 \left(\int_0^s |r-\tau|^{-2\alpha} d\tau\right)^{1/2}\le  \\
&
\|f(t;\cdot, s)\|_\infty  +  \|f(t;\cdot, s)\|_\infty   \sqrt{t} \frac 1 {1-2\alpha} \big(r^{1-2\alpha} +|s-r|^{1-2\alpha}\big)^{1/2} \le \\
&
\Big(1+\frac {\sqrt 2} {1-2\alpha}  t^{1-\alpha}\Big)\|f(t;\cdot, s)\|_\infty,
\end{align*}
which verifies the claimed bound.

\medskip
\noindent
\eqref{ii} The difference $\widetilde Q(r,s) := Q(r,s)-f(t;r,s)$ solves the equation 
\begin{equation}\label{Qteq}
\widetilde Q(r,s) + b_\alpha \int_0^s \widetilde Q(\tau,s)|r-\tau|^{-\alpha} d\tau = - b_\alpha \int_0^s |r-\tau|^{-\alpha} f(t; \tau,s)d\tau.
\end{equation}
For $r<s$, the expression in the right hand side  admits  the bound 
\begin{align*}
&
\left|
 b_\alpha \int_0^s |r-\tau|^{-\alpha} f(t; \tau,s)d\tau  \right| \le \\
&
  C(s,t)   \int_0^s |\tau-r|^{-\alpha} (s-\tau)^{-\alpha} d\tau =
  C(s,t)  \int_0^s v^{-\alpha} |s-r-v|^{-\alpha} dv = \\
&  C(s,t)  (s-r)^{1-2\alpha} \int_0^{s/(s-r)}\tau^{-\alpha}|1-\tau|^{-\alpha} d\tau =\\
&
  C(s,t)  (s-r)^{1-2\alpha}
\left(
\int_0^{s/(s-r)}\tau^{-\alpha}\big(|1-\tau|^{-\alpha}-\tau^{-\alpha}\big) d\tau 
+
\int_0^{s/(s-r)}\tau^{-2\alpha}   d\tau 
\right)\le \\
&
  C(s,t)  (s-r)^{1-2\alpha}
\left(
\int_0^{\infty}\tau^{-\alpha}\big||1-\tau|^{-\alpha}-\tau^{-\alpha}\big| d\tau 
+
\frac 1{1-2\alpha} \left(\frac{s}{s-r}\right)^{1-2\alpha}
\right) \le  \\
&
  C(s,t) \widetilde c(\alpha)
\Big( 
(s-r)^{1-2\alpha}+ s^{1-2\alpha}
\Big) \le     \widetilde c(\alpha)t^{1-2\alpha} C(s,t) =: C(t) C(s,t),
\end{align*}
where $\widetilde c(\alpha)$ is a positive constant.
The same bound holds for $r>s$. 
Thus, due to \eqref{i},
$$
|\widetilde Q(r,s)|\le C_1(t)C(t) C(s,t)
$$
and, consequently, 
\begin{align*}
|Q(r,s)| \le\, & |\widetilde Q(r,s)|+|f(t;r,s)| \le C_1(t)C(t) C(s,t)+ C(s,t) |s-r|^{-\alpha} \le  \\
&  
\big(
C_1(t)C(t) t^\alpha   + 1   
\big)C(s,t)|s -r |^{-\alpha} =: C_2(t) C(s,t) |s -r |^{-\alpha},
\end{align*}
 as claimed in \eqref{Qstalpha}.
  
\medskip
\noindent 
\eqref{iii}. 
In this case, the right hand side of \eqref{Qteq} satisfies 
\begin{align*}
&
\left|
 b_\alpha \int_0^s |r-\tau|^{-\alpha} f(t; \tau,s)d\tau  \right| \le\\
&
  \int_0^s |r-\tau|^{-\alpha} \big((s-\tau)^{-\alpha}-(t-\tau)^{-\alpha}\big)d\tau  =\\
&
\int_0^s |s-r-u|^{-\alpha} \big(u^{-\alpha} -(t-s+u)^{-\alpha}\big)du =\\
&
(t-s)^{1-2\alpha}\int_0^{s/(t-s)}  \Big|\frac{s-r}{t-s}-v\Big|^{-\alpha} \big(v^{-\alpha}-(v+1)^{-\alpha}\big)dv \le \\
&
(t-s)^{1-2\alpha}\int_0^{\infty}  |A-v|^{-\alpha} \big(v^{-\alpha}-(v+1)^{-\alpha}\big)dv,
\end{align*}
where we defined $A := (s-r)/(t-s)\in (0,\infty)$. 
 
Partition the integration region: 
\begin{equation}\label{split}
\begin{aligned}
&
\int_0^{\infty}  |A-v|^{-\alpha} \big(v^{-\alpha}-(v+1)^{-\alpha}\big)dv = \\
&
\int_0^{A}  (A-v)^{-\alpha} \big(v^{-\alpha}-(v+1)^{-\alpha}\big)dv
+
\int_A^{\infty}  (v-A)^{-\alpha} \big(v^{-\alpha}-(v+1)^{-\alpha}\big)dv.
\end{aligned}
\end{equation}
For any $A>0$, the first integral satisfies:
\begin{equation}\label{Aeq}
\begin{aligned}
&
\int_0^{A}  (A-v)^{-\alpha} \big(v^{-\alpha}-(v+1)^{-\alpha}\big)dv = \\
&
A^{1-2\alpha} \int_0^1 (1-u)^{-\alpha} \big(u^{-\alpha}-(u+1/A)^{-\alpha}\big)du \le 
  \int_0^1 (1-u)^{-\alpha} u^{ \alpha-1}du,
\end{aligned}
\end{equation}
where we used the bound
$$
x^{-\alpha}-y^{-\alpha} \le x^{-\alpha-\beta} (y-x)^\beta, \quad \forall y>x>0,\ \beta \in [0,1],\ \alpha \in [0,1),
$$
with $\beta :=1-2\alpha$. 
The second integral in \eqref{split} is also bounded: 
\begin{align*}
&
\int_A^{\infty}  (v-A)^{-\alpha} \big(v^{-\alpha}-(v+1)^{-\alpha}\big)dv =\\
&
 \int_A^{\infty}  (v-A)^{-\alpha}  \alpha \int_0^1(v+\tau)^{-\alpha -1}d\tau dv  =\\
&
\alpha \int_0^1 \int_0^{\infty}  u^{-\alpha}  (u+A+\tau)^{-\alpha -1} du d\tau = \\
&
\alpha \int_0^1 (A+\tau)^{-2\alpha }d\tau\int_0^{\infty}  w^{-\alpha}  (w +1)^{-\alpha -1} dw   \le  \\
&
\alpha  B(1-\alpha, 2\alpha) \int_0^1  \tau^{-2\alpha }d\tau = \frac{\alpha}{1-2\alpha}B(1-\alpha, 2\alpha).
\end{align*}
It follows that 
$$
\left|
 b_\alpha \int_0^s |r-\tau|^{-\alpha} f(t; \tau,s)d\tau  \right| \le C_\alpha (t-s)^{1-2\alpha}
$$
where 
$$
C_\alpha := \sup_{A>0}\int_0^{\infty}  |A-v|^{-\alpha} \big(v^{-\alpha}-(1+v)^{-\alpha}\big)dv<\infty.
$$
Applying \eqref{i} of Lemma \ref{lem} to equation \eqref{Qteq} yields the bound 
$$
|\widetilde Q(r,s)| \le C_1(t) C_\alpha (t-s)^{1-2\alpha},
$$
and, in turn,
$$
|Q(r,s)| \le\,    |\widetilde Q(r,s)| + |f(t;r,s)| \le 
 C_3(t) \big((t-s)^{1-2\alpha} +(s -r)^{-\alpha}- (t -r)^{-\alpha}\big)
$$
with $C_3(t) := 1+C_1(t) C_\alpha$.

\end{proof}
 
\subsection{Proof of \eqref{proveme}}
For $s<t$,
\begin{equation}\label{phisq}
\begin{aligned}
&
\E (\varphi_t(X)-\varphi_s(X))^2= \E\varphi_t(X)^2 + \E\varphi_s(X)^2 - 2\E \varphi_s(X)\varphi_t(X) =\\
&
\ \phantom +\int_0^t L(\tau,t) \left(L(\tau,t)+b_\alpha \int_0^t L(r,t)|\tau-r|^{-\alpha}dr\right)d\tau  
\\
&
+\int_0^s L(\tau,s) \left(L(\tau,s)+b_\alpha \int_0^s L(r,s) |r-\tau|^{-\alpha}dr\right)d\tau
\\
&
-2 \int_0^s L(\tau,s) \left(L(\tau,t)+b_\alpha \int_0^t L(r,t)|r-\tau|^{-\alpha} dr\right) d\tau =\\
&
- \int_0^t L(\tau,t) b_\alpha |\tau-t|^{-\alpha} d\tau  
-\int_0^s L(\tau,s) b_\alpha |\tau-s|^{-\alpha}d\tau \\
&
+2 \int_0^s L(\tau,s) b_\alpha |\tau-t|^{-\alpha}d\tau =: - b_\alpha (I_1+I_2+I_3),
\end{aligned}
\end{equation}
where we defined 
\begin{equation}\label{threeIs}
\begin{aligned}
I_1 := & \int_0^s \big(L(\tau,t)-L(\tau,s)\big)(t-\tau)^{-\alpha}d\tau, \\
I_2 := & \int_0^s L(\tau,s) \big((s-\tau)^{-\alpha}-(t-\tau)^{-\alpha}\big)d\tau, \\
I_3 := & \int_s^t L(\tau,t)(t-\tau)^{-\alpha} d\tau.
\end{aligned}
\end{equation}
We will estimate these terms using the bounds from Lemma \ref{lem}. By \eqref{ii} the solution to \eqref{Leq}
satisfies $|L(r,s)|\le  C_2(t) b_\alpha |r-s|^{-\alpha}$ and hence 
\begin{equation}\label{I3bnd}
|I_3| \le   C_2(t) \int_s^t  (t-\tau)^{-2\alpha}d\tau  = 
C_2(t) \frac 1 {1-2\alpha}(t-s)^{1-2\alpha}.
\end{equation}
Similarly,
\begin{equation}\label{I2bnd}
\begin{aligned}
|I_2| \le\ &   C_2(t)\int_0^s  (s-\tau)^{-\alpha}  \big((s-\tau)^{-\alpha}-(t-\tau)^{-\alpha}\big)d\tau =\\
&
C_2(t)\int_0^s  u^{-\alpha} \big(u^{-\alpha}-(t-s+u)^{-\alpha}\big)d\tau \le \\
&
C_2(t)(t-s)^{1-2\alpha}\int_0^\infty  v^{-\alpha} \big(v^{-\alpha}-(v+1)^{-\alpha}\big)dv= \\
&
C_2(t) \frac \alpha {1-2\alpha}B(1-\alpha, 2\alpha)(t-s)^{1-2\alpha}.
\end{aligned}
\end{equation}
To estimate $I_1$, note that $D(r,s):=L(r,t)-L(r,s)$ solves the equation, cf. \eqref{Leq}, 
\begin{multline}\label{dLeq}
D(r,s) + \int_0^s D(\tau,s) b_\alpha|r-\tau|^{-\alpha}d\tau  
= \\  b_\alpha \big((s-r)^{-\alpha}- (t-r)^{-\alpha}\big) -  \int_s^t L(\tau,t) b_\alpha|r-\tau|^{-\alpha}d\tau.
\end{multline}
By \eqref{ii} of Lemma \ref{lem},
\begin{align*}
&
\left|\int_s^t L(\tau,t) b_\alpha|r-\tau|^{-\alpha}d\tau\right|\le 
C_2(t) \int_s^t       (t-\tau)^{-\alpha} |r-\tau|^{-\alpha}d\tau =\\
&
C_2(t)\int_0^{t-s} u^{-\alpha} |t-r-u|^{-\alpha}du = 
C_2(t) (t-s)^{1-2\alpha} \int_0^1 v^{-\alpha} \left(\frac{t-r}{t-s}-v\right)^{-\alpha}dv \stackrel \dagger \le \\
&
C_2(t) (t-s)^{1-2\alpha} \left(\frac{t-r}{t-s}-1\right)^{-\alpha} \int_0^1 v^{-\alpha} dv =
C_2(t) \frac 1{1-\alpha} (t-s)^{1- \alpha}  (s-r)^{-\alpha},
\end{align*}
where $\dagger$ holds since $(t-r)/(t-s)\in [1,\infty)$ for $r\le s$. 
By linearity of \eqref{dLeq} and uniqueness of its solution, it follows from \eqref{ii} and \eqref{iii} of Lemma \ref{lem} that  
$$
|L(r,t)-L(r,s)|\le C_3(t) \big((t-s)^{1-2\alpha}+(s-r)^{-\alpha} -(t-r)^{-\alpha}\big) +  C_2(t)^2 \frac 1{1-\alpha} (t-s)^{1- \alpha}  (s-r)^{-\alpha}.
$$
Substitution of this bound yields
\begin{equation}\label{I1prebnd}  
\begin{aligned}
|I_1| \le\ & C_3(t) (t-s)^{1-2\alpha}\frac{t^{1-\alpha}}{1-\alpha} + C_3(t) \int_0^s  \big((s-\tau)^{-\alpha} -(t-\tau)^{-\alpha}\big)(t-\tau)^{-\alpha}d\tau 
+ \\
&
 C_2(t)^2 \frac 1{1-\alpha}(t-s)^{1-\alpha} \int_0^s  (s-\tau)^{-\alpha}(t-\tau)^{-\alpha}d\tau.
\end{aligned}
\end{equation}
The first integral in the right hand side satisfies the bound  
\begin{align*}
&
\int_0^s  \big((s-\tau)^{-\alpha} -(t-\tau)^{-\alpha}\big)(t-\tau)^{-\alpha}d\tau \le \\
&
\int_0^s  \big(u^{-\alpha} -(t-s+u)^{-\alpha}\big)(t-s+u)^{-\alpha}du \le\\
&
(t-s)^{1-2\alpha}\int_0^{\infty}  \big(v^{-\alpha} -(v+1)^{-\alpha}\big)(v+1)^{-\alpha}dv \le \\
&
\frac 1{1-2\alpha}(t-s)^{1-2\alpha}.
\end{align*}
The second integral in \eqref{I1prebnd} admits a similar estimate:
\begin{align*}
&
(t-s)^{1-\alpha} \int_0^s  (s-\tau)^{-\alpha}(t-\tau)^{-\alpha}d\tau \le \\
&
(t-s)^{1-\alpha} \int_0^s  (s-\tau)^{-2\alpha}d\tau = \frac 1{1-2\alpha} (t-s)^{1-\alpha}s^{1-2\alpha } \le \\
&
\frac 1{1-2\alpha} t^{1-\alpha } (t-s)^{1-2\alpha}.
\end{align*}
Thus we obtain 
\begin{equation}\label{I1bnd}
|I_1|\le 
\left(C_3(t) \frac{t^{1-\alpha}}{1-\alpha} +
C_3(t)  \frac 1{1-2\alpha}
+  
 C_2(t)^2  \frac 1{1-2\alpha} \frac 1{1-\alpha} t^{1- \alpha }\right)  (t-s)^{1-2\alpha}.
\end{equation}
Substitution of \eqref{I3bnd}, \eqref{I2bnd} and \eqref{I1bnd} into \eqref{phisq} yields \eqref{proveme}.

\section{Proof of \eqref{lower}}
Define the kernel
\begin{equation}\label{tildeL}
 \widetilde L(r,s) := L(r,s) + b_\alpha |s-r|^{-\alpha},
\end{equation}
where $L(r,s)$ solves equation \eqref{Leq}. 
Then the process \eqref{phi} can be written as $\varphi_s = \psi_s +\widetilde\varphi_s$, where 
$$
\psi_s = -\int_0^s b_\alpha |s-r|^{-\alpha}dX_r
\quad\text{and}\quad 
\widetilde\varphi_s:=\int_0^s \widetilde L(r,s) dX_r.
$$
To verify \eqref{lower} it suffices to show that, for some $t>0$,
\begin{equation}\label{lowb}
\varliminf_{s\nearrow t}\frac{\E\big(\psi_t  -\psi_s \big)^2}{(t-s)^{1-2\alpha}}  >0
\end{equation}
and 
\begin{equation}\label{upb}
\lim_{s\nearrow t}\frac{\E \big(\widetilde \varphi_t-\widetilde \varphi_s\big)^2}{(t-s)^{1-2\alpha}} =0.
\end{equation}

\subsection{Proof of \eqref{lowb}}
Denote $K(s,t) := \E  \psi_s\psi_t$ and write 
\begin{equation}\label{psist}
\begin{aligned}
\E (\psi_t-\psi_s)^2 =\ &  K(t,t)+K(s,s)-2K(s,t) = \\
& \big(K(t,t)-K(s,s)\big)-2\big(K(s,t)-K(s,s)\big).
\end{aligned}
\end{equation}
For $s<t$,  
\begin{align*}
K(s,t)  =\, & 
b_\alpha^2 \int_0^s   |s-r|^{-\alpha} |t-r|^{-\alpha}dr
+
b_\alpha^3 \int_0^s  \int_0^t |s-r|^{-\alpha}|t-\tau|^{-\alpha}|r-\tau|^{-\alpha} d\tau dr=\\ 
&
b_\alpha^2 \int_0^s   u^{-\alpha} |t -s+u|^{-\alpha}du
+
b_\alpha^3 \int_0^s \int_0^t u^{-\alpha}v^{-\alpha}|t-s+u -v|^{-\alpha}dvdu .
\end{align*}
In particular, 
\begin{align*}
K(t,t) =\, & 
b_\alpha^2 \int_0^t    u^{-2\alpha}  du
+
b_\alpha^3 \int_0^t u^{-\alpha} \int_0^t v^{-\alpha}|u -v|^{-\alpha}du dv =\\
&
\frac{b_\alpha^2}{1-2\alpha}      t^{1-2\alpha}  
+
\frac{2b_\alpha^3 B(1-\alpha,1-\alpha)}{2-3\alpha}t^{2-3\alpha} =: k_1 t^{1-2\alpha} +k_2 t^{2-3\alpha},
\end{align*}
and, consequently, for any $t>0$, 
$$
\lim_{s\nearrow t}{(t-s)^{1-2\alpha }}\big(K(t,t)-K(s,s)\big) = 0. 
$$
The second term in the right hand side of \eqref{psist} satisfies 
\begin{equation}\label{three}
\begin{aligned}
K(s,t)-K(s,s) = \ &
b_\alpha^2 \int_0^s    u^{-\alpha} \big( (t-s+u)^{-\alpha}-u^{-\alpha}\big)du +  \\
&
b_\alpha^3 \int_0^s \int_s^tu^{-\alpha}  v^{-\alpha}|t-s+u -v|^{-\alpha} dvdu + \\
&
b_\alpha^3 \int_0^s \int_0^su^{-\alpha}  v^{-\alpha}\big(|t-s+u -v|^{-\alpha}-| u -v|^{-\alpha}\big) dv du.
\end{aligned} 
\end{equation}
The first integral satisfies 
\begin{align*}
&
\frac 1{(t-s)^{1-2\alpha}}\int_0^s    u^{-\alpha}\big( (t-s+u)^{-\alpha}-u^{-\alpha}\big) du = \\
&
 \int_0^{s/(t-s)}    \tau^{-\alpha}\big( (1+\tau)^{-\alpha}- \tau^{-\alpha}\big) d \tau \xrightarrow[s\nearrow t]{} \\
 &
 \int_0^\infty   \tau^{-\alpha}\big( (1+\tau)^{-\alpha}- \tau^{-\alpha}\big) d \tau = -\frac \alpha {1-2\alpha}B(1-\alpha, 2\alpha).
\end{align*}
For $t>0$ and all sufficiently small $t-s>0$, the second integral can be written as   
\begin{align*}
& 
\int_0^s \int_s^tu^{-\alpha}  v^{-\alpha}|t-s+u -v|^{-\alpha} dvdu \le \\
&
s^{-\alpha} \int_0^s (s-r)^{-\alpha} \int_0^{t-s} | \tau -r  |^{-\alpha} d\tau dr =\\
& 
s^{-\alpha} \int_0^{t-s} (s-r)^{-\alpha} \int_0^{t-s} | \tau -r  |^{-\alpha} d\tau dr 
+\\
&
s^{-\alpha} \int_{t-s}^s (s-r)^{-\alpha} \int_0^{t-s} (r-\tau)^{-\alpha} d\tau dr =: J_1 +J_2,
\end{align*}
where   
$$
J_1 \le\,  
s^{-\alpha} (2s-t)^{-\alpha} \int_0^{t-s}  \int_0^{t-s} | \tau -r  |^{-\alpha} d\tau dr \le 
2 s^{-\alpha} (2s-t)^{-\alpha}    \frac 1 {1-\alpha}   (t-s)^{2-\alpha},
$$
and  
\begin{align*}
J_2 =\, & s^{-\alpha} \int_{t-s}^s (s-r)^{-\alpha} 
\frac 1{1-\alpha} \big(r^{1-\alpha}-(r-(t-s))^{1-\alpha}\big)dr \le \\
&
\frac 1{1-\alpha} s^{-\alpha}(t-s)^{1-\alpha} \int_{t-s}^s (s-r)^{-\alpha}  dr \le 
\frac 1{(1-\alpha)^2} s^{-\alpha}(t-s)^{1-\alpha} (2s-t)^{1-\alpha}.
\end{align*}
 Consequently, the second integral in \eqref{three} is of order $o((t-s)^{1-2\alpha})$ as $t-s\to 0$ and is therefore negligible.

To estimate the contribution of the third integral in \eqref{three}, let us partition its  domain into two regions so that
\begin{align*}
& 
\int_0^s \int_0^su^{-\alpha}  v^{-\alpha}\big(|t-s+u -v|^{-\alpha}-| u -v|^{-\alpha}\big) dv du = \\
&
\int_0^s u^{-\alpha} \int_0^u    v^{-\alpha}\big((t-s+u -v)^{-\alpha}-( u -v)^{-\alpha}\big) dv du
+\\
&
\int_0^s v^{-\alpha} \int_0^v  u^{-\alpha}  \big(|t-s+u -v|^{-\alpha}-(v -u)^{-\alpha}\big)  dudv =: I_1+I_2.
\end{align*}
Changing the integration variables in the first integral we get 
\begin{align*}
I_1 = &
\int_0^s u^{-\alpha} \int_0^u    (u-w)^{-\alpha}\big((t-s+w)^{-\alpha}-w^{-\alpha}\big) dw du = \\
&
(t-s)^{2-3\alpha }\int_0^{s/(t-s)}r^{-\alpha} \int_0^{r}    (r -\tau )^{-\alpha}\big((1+\tau)^{-\alpha}-\tau^{-\alpha}\big) d\tau  dr.
\end{align*}
The function 
$$
h(r) := \int_0^{r}    (r -\tau )^{-\alpha}\big(\tau^{-\alpha}-(\tau+1)^{-\alpha}\big) d\tau
$$
is continuous, non-negative and, as shown in \eqref{Aeq}, bounded. Since the function $f(\tau):= \tau^{-\alpha}-(\tau+1)^{-\alpha}$
is decreasing, integrable on $\Real$ and $f(\tau)\le \alpha\tau^{-\alpha-1}$, for any $r>0$
\begin{align*}
h(r) = &
\int_0^{r/2}    (r -\tau )^{-\alpha}f(\tau) d\tau
+
\int_{r/2}^r    (r -\tau )^{-\alpha}f(\tau) d\tau \le \\
&
(r/2 )^{-\alpha} \int_0^{\infty}    f(\tau) d\tau
+
f(r/2)\int_{r/2}^r    (r -\tau )^{-\alpha}  d\tau \le \\
&
(r/2 )^{-\alpha} \int_0^{\infty}    f(\tau) d\tau
+
\frac{\alpha}{1-\alpha}     (r/2)^{ -2\alpha} \le c(\alpha) r^{-\alpha}.
\end{align*}
Consequently 
$$
\frac{1}{(t-s)^{1-2\alpha}}|I_1| \le   
(t-s)^{1-\alpha }c(\alpha)\int_0^{s/(t-s)}r^{-2\alpha}    dr =(t-s)^{\alpha }\frac{c(\alpha)}{1-2\alpha} s^{1-2\alpha} 
 \xrightarrow[s\nearrow t]{}0.
$$
Similar bound holds for $I_2$ and therefore the third integral in \eqref{three} is negligible as well. 
Substituting all the above limits into \eqref{psist} we obtain:
$$
\frac {\E (\psi_t-\psi_s)^2}{(t-s)^{1-2\alpha}}
\xrightarrow[s\nearrow t]{}\      \frac {2b_\alpha^2\alpha} {1-2\alpha}B(1-\alpha, 2\alpha) >0,
$$
which verifies \eqref{lowb}.

\subsection{Proof of \eqref{upb}}
Substitution of \eqref{tildeL} into equation \eqref{Leq} shows that $\widetilde L(r,s)$ solves the integral equation
\begin{equation}\label{tildeLeq}
\widetilde L(r,s) +  \int_0^s  \widetilde L(\tau,s)  b_\alpha|r-\tau|^{-\alpha}d\tau 
= \Psi(r,s), \quad 0<r<s.
\end{equation}
where we defined 
\begin{equation}\label{Psidef}
\Psi(r,s):= b_\alpha^2 \int_0^s  (s-\tau)^{-\alpha} |r-\tau|^{-\alpha}d\tau.
\end{equation}
Note that this equation differs from \eqref{Leq} by its right hand side which in this case, being an iteration of a weakly singular kernel, is a locally bounded function:
\begin{equation}\label{Psibnd}
\begin{aligned}
&
\Psi(r,s) \le  \int_0^s  u^{-\alpha} |s-r-u|^{-\alpha}du = \\
&
\int_0^{s-r}  u^{-\alpha} (s-r-u)^{-\alpha}du
+
\int_{s-r}^s  u^{-\alpha} (u-(s-r))^{-\alpha}du \le \\
&
(s-r)^{1-2\alpha}B(1-\alpha, 1-\alpha)  
+
\int_0^r  (v+s-r)^{-\alpha} v^{-\alpha}dv \le \\
&
(s-r)^{1-2\alpha}B(1-\alpha, 1-\alpha)  
+
\frac 1{1-2\alpha} r^{1-2\alpha}\le C s^{1-2\alpha}.
\end{aligned}
\end{equation}
Similarly to \eqref{phisq},
$$
\E (\widetilde \varphi_t(X)-\widetilde \varphi_s(X))^2= I_3(s,t) +  I_1(s,t)   +I_2(s,t),
$$
where,  cf. \eqref{threeIs}, 
\begin{align*}
I_1(s,t) := & \int_0^s \big(\widetilde L(\tau,t)-\widetilde L(\tau,s)\big)\Psi(\tau,t) d\tau , \\
I_2(s,t) := & \int_0^s \widetilde L(\tau,s) \big(\Psi(\tau,s) - \Psi(\tau,t)\big) d\tau, \\
I_3(s,t) := & \int_s^t \widetilde L(\tau,t) \Psi(\tau,t) d\tau.
\end{align*}
It remains to check that each one of these integrals are of order $o\big((t-s)^{1-2\alpha}\big)$ as $s\to t$.  
By definition \eqref{tildeL}, the function $\widetilde L(r,s)$ is bounded:
\begin{equation}\label{tildeLbnd}
\begin{aligned}
&
|\widetilde L(r,s)|  =  \left|\int_0^s L(\tau,s) b_\alpha |r-\tau|^{-\alpha}d\tau\right|  \le 
\|L(\cdot,s)\|_2 \left(\int_0^s   |r-\tau|^{-2\alpha}d\tau\right)^{1/2}= \\
&
\|L(\cdot,s)\|_2 \frac 1{\sqrt{1-2\alpha}} \left( r^{1-2\alpha}+(s-r)^{1-2\alpha}\right)^{1/2} \le
\|L(\cdot,s)\|_2 \frac 1{ \sqrt{1/2 -\alpha} }    s^{ 1/2- \alpha} =: A(s).
\end{aligned}
\end{equation}
Due to \eqref{Psibnd} and \eqref{tildeLbnd},
$$
 I_3(s,t)  =  O(t-s) = o\big((t-s)^{1-2\alpha}\big).
$$
In view of \eqref{tildeLbnd}, 
\begin{equation}\label{absI2}
|I_2(s,t)| \le\,   A(s) \int_0^s  \big|\Psi(r,s) - \Psi(r,t)\big| dr 
\end{equation}
where, by definition \eqref{Psidef},
\begin{equation}\label{I2}
\begin{aligned}
  &   \int_0^s  \big|\Psi(r,s) - \Psi(r,t)\big| dr \le \\
&
  \int_0^s  \left|\int_0^s  (s-\tau)^{-\alpha} |r-\tau|^{-\alpha}d\tau
 - \int_0^t  (t-\tau)^{-\alpha} |r-\tau|^{-\alpha}d\tau\right| dr \le \\
&
   \int_0^s   \int_0^s  \big((s-\tau)^{-\alpha} -(t-\tau)^{-\alpha}\big)|r-\tau|^{-\alpha}d\tau
   dr
 + \\
& 
   \int_0^s \int_s^t  (t-\tau)^{-\alpha} |r-\tau|^{-\alpha}d\tau dr =:  J_1+J_2 .
\end{aligned}
\end{equation}
The first term in the right hand side satisfies  
\begin{align*}
(1-\alpha) J_1 = &   (1-\alpha) \int_0^s  \big((s-\tau)^{-\alpha} -(t-\tau)^{-\alpha}\big)\Big(\int_0^s|r-\tau|^{-\alpha}dr\Big)d\tau
   =\\
&  
\int_0^s \Big( (s-\tau)^{-\alpha}-(t-\tau)^{-\alpha}\Big) \big(\tau^{1-\alpha}+(s-\tau)^{1-\alpha}\big)d\tau = \\
&
\int_0^s  (s-\tau)^{-\alpha}  \tau^{1-\alpha} d\tau
-
\int_0^t  (t-\tau)^{-\alpha}   \tau^{1-\alpha} d\tau
+
\int_s^t  (t-\tau)^{-\alpha}   \tau^{1-\alpha} d\tau
+\\
&
\int_0^s \Big( (s-\tau)^{-\alpha}-(t-\tau)^{-\alpha}\Big)  (s-\tau)^{1-\alpha} d\tau.
\end{align*} 
The second term can be written as  
\begin{align*}
(1-\alpha)J_2 =   (1-\alpha) & \int_s^t  (t-\tau)^{-\alpha} \Big(\int_0^s(\tau-r)^{-\alpha} dr \Big)d\tau =\\
&
   \int_s^t  (t-\tau)^{-\alpha} \big(\tau^{1-\alpha}-(\tau-s)^{1-\alpha}\big) d\tau.
\end{align*}
The integrals in these expressions satisfy 
\begin{multline*}
\int_0^t  (t-\tau)^{-\alpha}   \tau^{1-\alpha} d\tau
-
\int_0^s  (s-\tau)^{-\alpha}  \tau^{1-\alpha} d\tau = \\
B(1-\alpha, 2-\alpha)\big(t^{2-2\alpha}-s^{2-2\alpha}\big) = O( t-s ),
\end{multline*} 
$$
\shoveright{
\int_s^t  (t-\tau)^{-\alpha}   \tau^{1-\alpha} d\tau =   \int_0^{t-s}  v^{-\alpha}   (t-v)^{1-\alpha} dv \le  
  \frac {t^{1-\alpha}}{1-\alpha} (t-s)^{1-\alpha} =
O((t-s)^{1-\alpha}), 
}
$$
\begin{multline*}
\int_0^s \Big( (s-\tau)^{-\alpha}-(t-\tau)^{-\alpha}\Big)  (s-\tau)^{1-\alpha} d\tau = 
\int_0^s \Big( v^{-\alpha}-(t-s+v)^{-\alpha}\Big)  v^{1-\alpha} dv =  \\
 (t-s)^{2-2\alpha}\int_0^{s/(t-s)}
\Big( \tau^{-\alpha}-(\tau+1)^{-\alpha}\Big)  \tau^{1-\alpha} d\tau = O(t-s),
\end{multline*}
\begin{multline*}
\int_s^t  (t-\tau)^{-\alpha}  (\tau-s)^{1-\alpha}  d\tau = \int_0^{t-s} (t-s-v)^{-\alpha}v^{1-\alpha}dv = \\
B(1-\alpha, 2-\alpha)
(t-s)^{2-2\alpha} = O((t-s)^{2-2\alpha}).
\end{multline*}
Substitution of these estimates into \eqref{I2} and \eqref{absI2} yields
$$
I_2(s,t)  = O\big(t-s\big) = o\big( (t-s)^{1-2\alpha} \big).
$$
It remains to bound $I_3(s,t)$. To this end,  in view of \eqref{tildeLeq}, 
\begin{multline*} 
\widetilde L(r,t) - \widetilde L(r,s) +  \int_0^s  \big(\widetilde L(\tau,t)-\widetilde L(\tau,s)\big)  b_\alpha|r-\tau|^{-\alpha}d\tau = \\
-
\int_s^t  \widetilde L(\tau,t)  b_\alpha|r-\tau|^{-\alpha}d\tau + \Psi(r,t)-\Psi(r,s), \quad 0<r<s<t.
\end{multline*}
Consequently,
\begin{align*}
\int_0^s \big|\widetilde L(r,t) - \widetilde L(r,s)\big|dr \le &
\int_0^s  \big|\widetilde L(\tau,t)-\widetilde L(\tau,s)\big|  b_\alpha\int_0^s |r-\tau|^{-\alpha}dr d\tau 
+\\
&
\int_s^t  |\widetilde L(\tau,t)|  b_\alpha\int_0^s |r-\tau|^{-\alpha}dr d\tau + 
\int_0^s |\Psi(r,t)-\Psi(r,s)|dr,
\end{align*}
where 
$$
b_\alpha\int_0^s |r-\tau|^{-\alpha}dr   = (1-\tfrac \alpha 2) \Big(\tau^{1-\alpha}+(s-\tau)^{1-\alpha}\Big) \le 2s^{1-\alpha} \le 2t^{1-\alpha}.
$$
Hence for sufficiently small $t>0$, 
\begin{align*}
|I_1(s,t)| \le & \sup_{\tau\le s}\Psi(\tau,t)
\int_0^s \big|\widetilde L(r,t) - \widetilde L(r,s)\big|dr \le \\
&
\frac{1}{1-2t^{1-\alpha}}
\Big(
A(t) 2t^{1-\alpha} (t-s)+ \int_0^s |\Psi(r,t)-\Psi(r,s)|dr
\Big) = O(t-s),
\end{align*}
where the last estimate holds since $\Psi(\tau,t)$ is bounded and, as we already argued above, cf. \eqref{I2},
$$
\int_0^s |\Psi(r,t)-\Psi(r,s)| = O(t-s).
$$
This completes the proof of \eqref{upb}.



\begin{thebibliography}{10}

\bibitem{BN03}
Fabrice Baudoin and David Nualart, \emph{Equivalence of {V}olterra processes},
  Stochastic Process. Appl. \textbf{107} (2003), no.~2, 327--350. \MR{1999794
  (2004e:60064)}

\bibitem{BSV07}
Christian Bender, Tommi Sottinen, and Esko Valkeila, \emph{Fractional processes
  as models in stochastic finance}, Advanced mathematical methods for finance,
  Springer, Heidelberg, 2011, pp.~75--103. \MR{2792076}

\bibitem{CCK}
Chunhao Cai, Pavel Chigansky, and Marina Kleptsyna, \emph{Mixed {G}aussian
  processes: a filtering approach}, Ann. Probab. \textbf{44} (2016), no.~4,
  3032--3075. \MR{3531685}

\bibitem{Ch01}
Patrick Cheridito, \emph{Mixed fractional {B}rownian motion}, Bernoulli
  \textbf{7} (2001), no.~6, 913--934. \MR{1873835 (2002k:60163)}

\bibitem{Ch03}
\bysame, \emph{Arbitrage in fractional {B}rownian motion models}, Finance
  Stoch. \textbf{7} (2003), no.~4, 533--553. \MR{2014249 (2004m:60076)}

\bibitem{Ch03b}
\bysame, \emph{Representations of {G}aussian measures that are equivalent to
  {W}iener measure}, S\'eminaire de {P}robabilit\'es {XXXVII}, Lecture Notes in
  Math., vol. 1832, Springer, Berlin, 2003, pp.~81--89. \MR{2053042}

\bibitem{ChKM20}
P.~Chigansky, M.~Kleptsyna, and D.~Marushkevych, \emph{Mixed fractional
  {B}rownian motion: a spectral take}, J. Math. Anal. Appl. \textbf{482}
  (2020), no.~2, 123558, 23. \MR{4015693}

\bibitem{CDM2025}
Carsten~H. Chong, Thomas Delerue, and Fabian Mies, \emph{Rate-optimal
  estimation of mixed semimartingales}, Ann. Statist. \textbf{53} (2025),
  no.~1, 219--244. \MR{4865014}

\bibitem{DSSV24}
Josephine Dufitinema, Foad Shokrollahi, Tommi Sottinen, and Lauri Viitasaari,
  \emph{Long-range dependent completely correlated mixed fractional brownian
  motion}, Stochastic Processes and their Applications \textbf{170} (2024),
  104289.

\bibitem{EM02}
Paul Embrechts and Makoto Maejima, \emph{Selfsimilar processes}, Princeton
  Series in Applied Mathematics, Princeton University Press, Princeton, NJ,
  2002. \MR{1920153 (2004c:60003)}

\bibitem{gatheral2022volatility}
Jim Gatheral, Thibault Jaisson, and Mathieu Rosenbaum, \emph{Volatility is
  rough}, Quantitative Finance \textbf{18} (2018), no.~6, 933--1336.

\bibitem{HH76}
Takeyuki Hida and Masuyuki Hitsuda, \emph{Gaussian processes}, Translations of
  Mathematical Monographs, vol. 120, American Mathematical Society, Providence,
  RI, 1993, Translated from the 1976 Japanese original by the authors.
  \MR{1216518 (95j:60057)}

\bibitem{Hi68}
Masuyuki Hitsuda, \emph{Representation of {G}aussian processes equivalent to
  {W}iener process}, Osaka J. Math. \textbf{5} (1968), 299--312. \MR{0243614
  (39 \#4935)}

\bibitem{jusselin2020no}
Paul Jusselin and Mathieu Rosenbaum, \emph{No-arbitrage implies power-law
  market impact and rough volatility}, Mathematical Finance \textbf{30} (2020),
  no.~4, 1309--1336.

\bibitem{M08}
Yuliya~S. Mishura, \emph{Stochastic calculus for fractional {B}rownian motion
  and related processes}, Lecture Notes in Mathematics, vol. 1929,
  Springer-Verlag, Berlin, 2008. \MR{2378138 (2008m:60064)}

\bibitem{muhle2025metaflow}
Johannes Muhle-Karbe, Youssef~Ouazzani Chahdi, Mathieu Rosenbaum, and Grégoire
  Szymanski, \emph{A unified theory of order flow, market impact, and
  volatility}, 2026, arXiv:2601.23172.

\bibitem{PT17}
Vladas Pipiras and Murad~S. Taqqu, \emph{Long-range dependence and
  self-similarity}, Cambridge Series in Statistical and Probabilistic
  Mathematics, [45], Cambridge University Press, Cambridge, 2017. \MR{3729426}

\bibitem{RN55}
Frigyes Riesz and B{\'e}la Sz.-Nagy, \emph{Functional analysis}, Dover Books on
  Advanced Mathematics, Dover Publications Inc., New York, 1990, Translated
  from the second French edition by Leo F. Boron, Reprint of the 1955 original.
  \MR{1068530 (91g:00002)}

\bibitem{Sh66}
L.~A. Shepp, \emph{Radon-{N}ikod\'ym derivatives of {G}aussian measures}, Ann.
  Math. Statist. \textbf{37} (1966), 321--354. \MR{0190999 (32 \#8408)}

\bibitem{vZ07}
Harry van Zanten, \emph{When is a linear combination of independent f{B}m's
  equivalent to a single f{B}m?}, Stochastic Process. Appl. \textbf{117}
  (2007), no.~1, 57--70. \MR{2287103 (2008h:60140)}

\end{thebibliography}

\def\cprime{$'$} \def\cprime{$'$} \def\cydot{\leavevmode\raise.4ex\hbox{.}}
  \def\cprime{$'$} \def\cprime{$'$} \def\cprime{$'$}
\providecommand{\bysame}{\leavevmode\hbox to3em{\hrulefill}\thinspace}
\providecommand{\MR}{\relax\ifhmode\unskip\space\fi MR }
\providecommand{\MRhref}[2]{%
  \href{http://www.ams.org/mathscinet-getitem?mr=#1}{#2}
}
\providecommand{\href}[2]{#2}

\begin{acks}
We are grateful to Mathieu Rosenbaum for bringing the question to our attention and for his insightful comments on the subject.
\end{acks}


\end{document}